\documentclass[journal]{IEEEtran}
\usepackage{cite}
\usepackage{amssymb}
\usepackage{amsthm}
\usepackage{stmaryrd}
\usepackage[nosumlimits]{amsmath}
\usepackage{algorithm}
\usepackage[noend]{algpseudocode}
\usepackage{graphicx}
\usepackage{textcomp}
\usepackage{setspace}
\usepackage{xcolor}
\usepackage{steinmetz}
\usepackage{hyperref} 
\usepackage{cleveref}
\usepackage{mathtools}
\usepackage{soul}

\usepackage{stackengine}

\stackMath

\usepackage{tikz}
\newcommand{\tikzmark}[1]{\tikz[overlay,remember picture] \node (#1) {};}
\tikzset{square arrow/.style={to path={-- ++(0,-.25) -| (\tikztotarget)}}}
\tikzset{deep square arrow/.style={to path={-- ++(0,-.5) -| (\tikztotarget)}}}

\DeclareMathOperator{\sign}{sign}

\usepackage{accents}

\usepackage{algorithm}

\title{Distribution System Flexibility Characterization:\\ A Network-Informed Data-Driven Approach}
\author{Qi~Li,~\IEEEmembership{Student~Member,~IEEE},~Jianzhe~Liu,~\IEEEmembership{Member,~IEEE},~Bai~Cui,~\IEEEmembership{Member,~IEEE},\\~Wenzhan~Song,~\IEEEmembership{Senior~Member,~IEEE},~and~Jin~Ye,~~\IEEEmembership{Senior~Member,~IEEE}
\thanks{Q. Li, W. Song, and J. Ye are with University of Georgia, Georgia, U.S.A.}
\thanks{J. Liu is with Shanghai Jiao Tong University, Shanghai, China.}
\thanks{B. Cui is with Iowa State University, Ames, IA, U.S.A..}
\thanks{J. Liu is the corresponding author.}
}
\allowdisplaybreaks[4]
\begin{document}

\maketitle
\begin{abstract}
    A distribution system can flexibly adjust its substation-level power output by aggregating its local distributed energy resources (DERs). Due to DER and network constraints, characterizing the exact feasible power output region is computationally intensive. 
    Hence, existing results usually rely on unpractical assumptions or suffer from conservativeness issues. 
    Sampling-based data-driven methods can potentially address these limitations. Still, existing works usually exhibit computational inefficiency issues as they use a random sampling approach, which carries little information from network physics and provides few insights into the iterative search process.
    This letter proposes a novel network-informed data-driven method to close this gap. 
    A computationally efficient data sampling approach is developed to obtain high-quality training data, leveraging network information and legacy learning experience. 
    Then, a classifier is trained to estimate the feasible power output region with high accuracy. 
    Numerical studies based on a real-world Southern California Edison network validate the performance of the proposed work. 
\end{abstract}

\section{Introduction}

Proper coordination of distributed energy resources (DERs) transforms a passive distribution system into an active grid asset. From the grid operation standpoint, it is critical to characterize the distribution system flexibility region, i.e., the set of feasible substation-level power outputs subject to network and component operational constraints. This set is essentially a projection from the high-dimensional DER and network operation region, and, in general, finding its exact characterization is computationally unrealistic~\cite{chen2021leveraging}.

A variety of approximation methods to characterize the flexibility region have been developed. 
For example, a Minkowski sum-based approximation method is proposed; this method is scalable but cannot handle network constraints~\cite{cui2020network}.
Ref.~\cite{chen2021leveraging,cui2020network} use robust optimization methods to find a convex inner approximation of the flexibility set, explicitly considering network constraints and the temporal coupling of the DER operation decisions. Nevertheless, these approximations are conservative, and the shapes of the approximated set are fixed and presumed, which do not necessarily correspond to the actual geometry. 
Data-driven methods have been investigated as well~\cite{ageeva2019analysis,taheri2022data}. They usually use a random sampling approach and numerical approaches based on iterative algorithms to find labeled data for training purposes. The sampling and labeling operations could limit the scalability and bring in high computational overhead.



\begin{figure}
\centering
\includegraphics[width=0.7\linewidth]{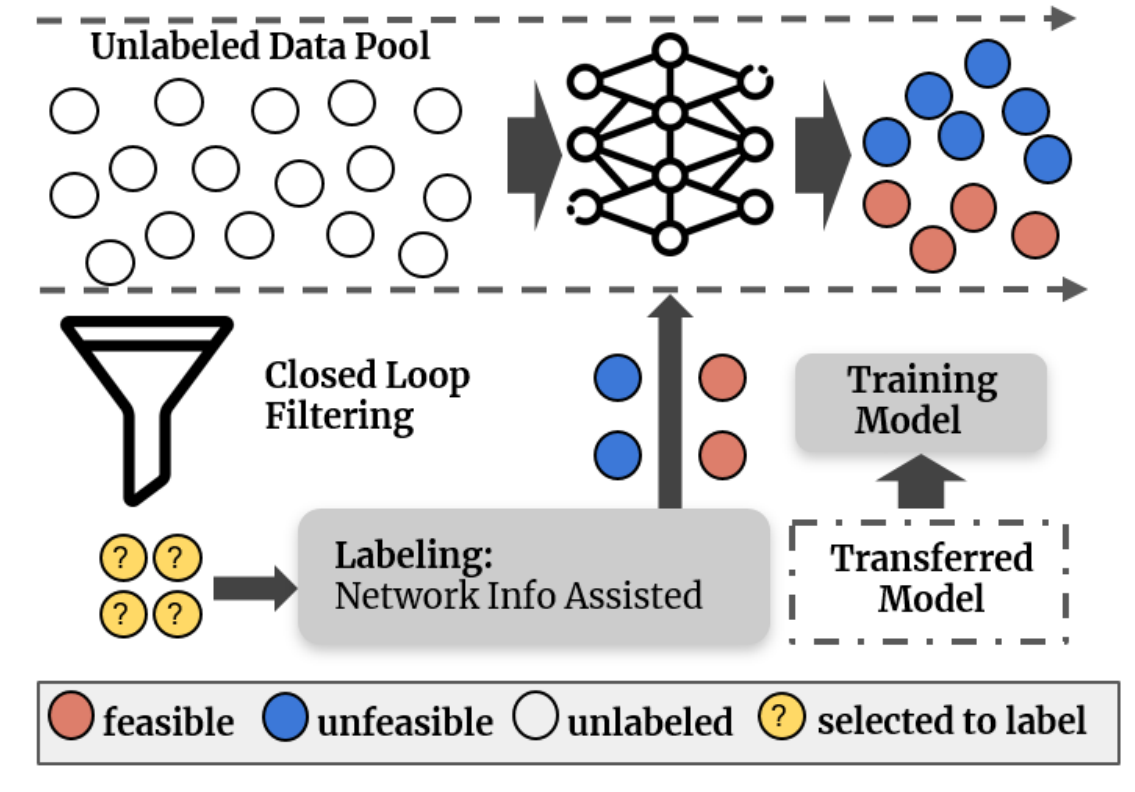}
\caption{Proposed training method to obtain the classifier}
\label{fig:architecture}
\vspace{-0.2in}
\end{figure}

To close this gap, we propose a network-informed data-driven approximation approach that exhibits superior scalability. Our main contributions are two-fold. First, unlike existing methods that use iterative algorithms or prescribed approximation shapes, we propose a new approach that uses a highly scalable matrix operation-based classifier to efficiently sketch an approximated region with limited conservativeness. Second, the classifier is obtained by a novel training strategy with high efficiency. As shown in Fig.~\ref{fig:architecture}, we develop a closed-loop data filtering algorithm to actively select samples that are most helpful to classifier training in a rolling horizon. Moreover, we explicitly use the network knowledge to develop a rigorous condition for sample labeling. This essentially trims the sample space to improve approximation accuracy and scalability further. 


\section{System Model and Problem Description}
We consider a distribution system with one substation feeder bus and $n$ load buses, on which there are $m$ controllable DERs. The time horizon is given by $\mathcal{T} = \{1,\cdots,T\}$. The power outputs of the DERs are managed and aggregated to achieve controllability for the substation-level power output so that the distribution system becomes a controllable grid asset.

\vspace{-0.10in}
\subsection{Distribution System Operation Model}\label{sec:model:model}
The DER aggregation considers 1) DER capacity limits represented by interval constraints; 2) network constraints, including linearized power flow equations and interval voltage limits. 
Based on the above discussion, the system operation constraints are modeled in the following compact form\cite{chen2021leveraging}:
\begin{align}
    &\mathbf{Wp \leq z}, \quad \mathbf{p_0 = Dp+b}, \label{eq:main}
\end{align}
where $\mathbf{p} \in \mathbb{R}^{mT}$ and $\mathbf{p_0} \in \mathbb{R}^{T}$ represent the dispatchable DER power outputs and substation-level power output, respectively; $\mathbf{W}$ and $\mathbf{D}$ are both given constant matrices such that $\mathbf{W}$ captures the DER operational constraints and the network voltage constraints, and $\mathbf{D}$ models the mapping of DER power outputs to the substation; $\mathbf{z}$ and~$\mathbf{b}$ are constant coefficients, representing given parameters such as load forecasts.

The inequality in~\eqref{eq:main} represents the DER and network operational constraints, and the equality constraint models the mapping of DER and load power to substation-level power output based on the linear power flow model. The use of the linear power flow model is justified owing to the tight voltage limits in the distribution system.
The constraint captures the steady-state behavior of various kinds of DERs, including HVACs, energy storage units, and photovoltaics~\cite{chen2021leveraging,cui2020network}. 

\vspace{-0.035in}
\subsection{Flexibility Characterization Problem}
Distribution system flexibility set (DSFS) refers to the set of all the substation-level power output realizations that are feasible to~\eqref{eq:main} with appropriate $\mathbf{p}$. Obtaining its exact characterization is generally computationally expensive. 
A network-informed data-driven approach is proposed in this letter. First, we use a novel offline training method to obtain a classifier that determines whether a substation-level power output sample belongs to the DSFS. Then, the samples from the substation-level power output space are classified; the union of the identified DSFS members forms a DSFS estimation. Note that the second step is scalable as it only involves a) computationally trivial sampling operations and b) simple matrix operations associated with the classification. Compared to an iterative algorithm-based numerical method, our method is five orders of magnitude faster, as shown in Section~\ref{sec:case}. 

Nonetheless, the offline training step to obtain the classifier is more computationally demanding. Developing a new efficient training strategy is the focus of this letter.

\vspace{-0.035in}
\section{Proposed Training Strategy}
Traditional data-driven classifier training strategy can be summarized as \emph{sampling $\to$ labeling $\to$ training}. It is an open-loop process where one needs to prepare a training dataset before training commences. Given no information about the sample space geometry, a larger set of randomly drawn samples is usually needed to ensure the representativeness of the sample space at the cost of increased computational burdens. 

To circumvent this issue, an active training strategy is proposed, as illustrated below:
\begin{equation}
    \text{sampling} \to \text{filt}\tikzmark{a}\text{ering} \to \text{lab}\tikzmark{b}\text{eling} \to \text{train}\tikzmark{c}\text{ing} \nonumber
    \tikz[overlay,remember picture]
   {\draw[->,square arrow] (c.south) to (a.south);}
   \tikz[overlay,remember picture]
   {\draw[->,square arrow] (c.south) to (b.south);}
\end{equation}

We create a closed-loop training process where the classifier is trained through multiple steps, as shown in Fig.~\ref{fig:architecture}: 1) sampling: randomly select samples from the unlabeled pool (colorless circles); 2) filtering: determine posterior probabilities and select the most uncertain samples (yellow circles with question marks) for labeling; 3) labeling: label selected samples as feasible (red circle) or unfeasible (blue circle), leveraging network knowledge; 4) training: train the model using the enlarged training set, including newly selected samples, in which transfer learning can be used to accelerate training, utilizing parameters from a historical model (represented by the dotted box in the figure).
Ideally, the dataset size is relatively small at first and then grows sequentially by incorporating selected high-value training data points identified in each epoch. Here we use the growing knowledge about the sample space to develop a filtering algorithm for such data selection. The filtering and labeling algorithms keep improving to ensure accuracy and scalability throughout the training process, as will be discussed later in detail. 
It is worth noting that although the feedback-learning framework is first proposed in the machine learning community~\cite{ren2021survey}, here, it is used as a vehicle to implement the nontrivial and novel network-informed algorithms.


\vspace{-0.15in}
\subsection{Network-Informed Labeling}
Each training data point consists of a substation-level power output sample and a label about whether this sample is feasible, i.e., belonging to DSFS. Let $\mathbf{x}_i=[\mathbf{\hat{p}_{0,i}}^\top,y_i]^\top$, where $\mathbf{\hat{p}_{0,i}}$ is the sample, and $y_i$ is the label with $1$ representing ``feasible'' and $0$ otherwise. In practice, this label is obtained through numerical methods to test whether a $\mathbf{\hat{p}_{0,i}}$ is feasible to~\eqref{eq:main}, which are computationally intensive when dealing with a large number of samples.

To simplify the process, we leverage the network knowledge to trim the sample space such that points from a certain region bear no need for numerical labeling. To this end, we first find a convex inner approximation of DSFS, whence any members must have a ``1'' label, by solving:
\begin{align} \label{eq:aro}
    \max_{\mathbf{p}^-_0 \le \mathbf{p}^+_0} \left\{ \mathbf{1}^\top \left( \mathbf{p}^+_0 - \mathbf{p}^-_0 \right) + \min_{\mathbf{p}^-_0 \le \mathbf{p}_0 \le \mathbf{p}^+_0} \max_{\substack{\mathbf{p}_0 = \mathbf{D}\mathbf{p} + \mathbf{b} \\ \mathbf{W}\mathbf{p} \le \mathbf{z}}} \mathbf{0}^\top \mathbf{p} \right\}
\end{align}
where $\mathbf{p}_0^+, \mathbf{p}_0^- \in \mathbb{R}^T$ represent the upper and lower bounds of the substation-level power output, respectively. The inner min-max (feasibility) problem admits the optimal value of $0$ if and only if for any substation-level power output between $\mathbf{p}_0^-$ and $\mathbf{p}_0^+$, there exists a DER output schedule that makes all the operational constraints described by~\eqref{eq:main} satisfied. Therefore, the hyperbox $\{\mathbf{p}_0: \mathbf{p}^-_0 \le \mathbf{p}_0 \le \mathbf{p}^+_0\}$ must be a subset of DSFS when the optimal value of the outer problem is finite. Note that \eqref{eq:aro} is an adaptive robust optimization (ARO) problem. 
One usually makes $\mathbf{p}$ a function of $\mathbf{p}_0$ in solving an ARO problem of this type. This paper assumes an affine decision rule. The problem then reformulates into a max-min problem in the form of $\max \min \mathbf{1}^\top \left( \mathbf{p}^+_0 - \mathbf{p}^-_0 \right)$. Inserting a slack variable $\mathbf{s} = \min \mathbf{1}^\top \left( \mathbf{p}^+_0 - \mathbf{p}^-_0 \right)$ yields a standard robust linear programming problem with the objective function becoming $\max \mathbf{s}$. The problem is tractable with well-established solution methods.


We enlarge the DSFS approximation and further trim the sample space in each epoch. Note that~\eqref{eq:main} are convex constraints, and DSFS is a projection of the feasibility region of~\eqref{eq:main} onto the $\mathbf{p}_0$-space; hence, DSFS is a convex set, and the convex hull of any DSFS members must be a subset of DSFS.
Recall that we train the classifier in epochs; in each epoch, the training set is expanded by new samples.
Given a DSFS subset in an epoch, we only need to label those lying outside of the subset numerically.
Then, the convex hull of those new samples with ``1'' labels and the original DSFS subset becomes a new DSFS subset. Hence, in the next epoch, those new samples lying in this enlarged set can be directly labeled again, thanks to the use of network information.





\subsection{Closed Loop Filtering}
\label{sec:closed-loop}
In each epoch, we seek to find the samples that are most \emph{uncertain} to the classifier, i.e., containing the most \emph{fresh} knowledge about the sample space. 

An uncertainty quantification method is applied. Let $\mathbb{P}(1|\mathbf{\hat{p}}_{0,i})$ be the posterior probability of a sample being feasible, according to an estimator. The closer $\mathbb{P}(1|\mathbf{\hat{p}}_{0,i})$ is to $1$ (resp., $0$), the more likely the sample is feasible (resp., infeasible); whereas the closer it is to $0.5$, the more uncertain it is. Then, by a simple mapping, we can find a monotone uncertainty metric: If $\mathbb{P}(1|\mathbf{\hat{p}}_{0,i}) > 0.5$, let $\mathbb{M}(\mathbf{\hat{p}}_{0,i}) = 2(1-\mathbb{P}(1|\mathbf{\hat{p}}_{0,i}))$; otherwise, $\mathbb{M}(\mathbf{\hat{p}}_{0,i}) = 2\mathbb{P}(1|\mathbf{\hat{p}}_{0,i})$, where $\mathbb{M}(\mathbf{\hat{p}}_{0,i})$ is the quantified uncertainty. After using this metric to evaluate all unlabeled samples, the most uncertain samples can be selected by ranking the quantified uncertainties.
As for the initial number of samples and the selected number of samples in each epoch, they are open to customization, which acts as the hyperparameters for our model. 

The execution of the aforementioned process depends on finding $\mathbb{P}(1|\mathbf{\hat{p}}_{0,i})$. The classifier is structured to accomplish this task. We build the classifier using a multi-layer perceptron (MLP) model, defined as $f(\mathbf{\hat{p}}_{0,i}): \mathbb{R}^T \to [0,1]$. Its output is $\mathbb{P}(1|\mathbf{\hat{p}}_{0,i})$. Meanwhile, if a classification result ($0$ or $1$) is needed, a simple probabilistic smoothing approximation can be used, for example, $\sign f$, which is $1$ if $f > 0.5$ and $0$ otherwise. It is worth mentioning that the proposed strategy is general, and we can use models other than the MLP model. 

With the above discussion, the closed-loop filtering is conducted as follows: In each epoch, given an unlabeled sample pool, we first find the posterior probability of each unlabeled sample using the classifier obtained in the last epoch (or the initial classifier); then, the most uncertain samples with a suitable size are selected to label and then train the classifier; the updated classifier is then similarly used in the next epoch. The initial classifier's parameters can either be randomly generated or transferred from a historical model. Numerical testing suggests that the transfer learning approach is effective in characterizing the DSFS, for the transferred model entails substation-level power output sample space geometry knowledge that can warm-start the training. 
In addition, the training speed per epoch is accelerated since there are fewer trainable parameters during transfer learning.

\section{Case Studies} \label{sec:case}
In this section, we conduct numerical testing based on a three-phase distribution feeder of Southern California Edison (SCE) with 126 load buses and 366 DER having temporal couplings~\cite{chen2021leveraging,cui2020network}. We estimate the aggregated flexibility region of the substation-level real power output profile. 

\begin{figure}[t]
\centering
\includegraphics[width=0.90\linewidth]{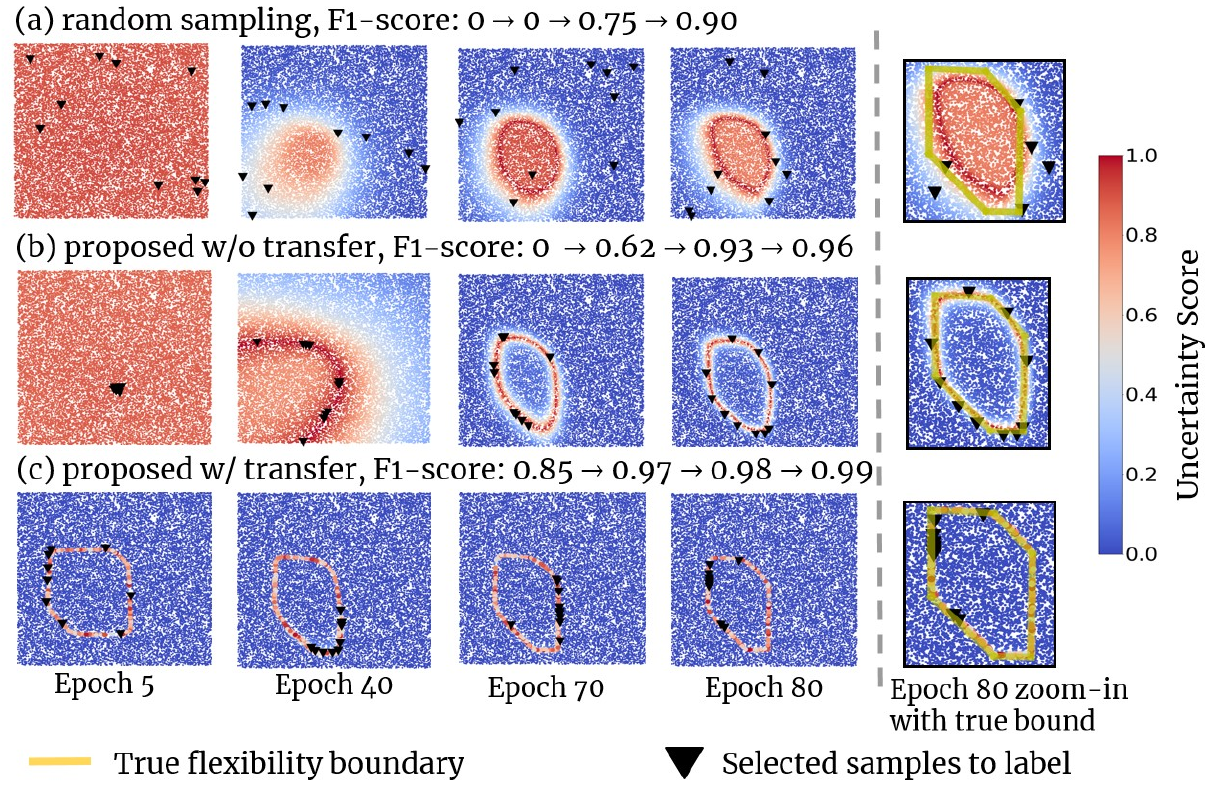}
\caption{Benchmarking the proposed method using uncertainty heatmap.}
\label{fig:case1}
\end{figure}

For visualization purposes, we first conduct a numerical study regarding a two-dimensional aggregated power profile. With a time step (TS) of one hour, the flexibility for the time window [8,10] is estimated. For the specific setting, we randomly picked 100 samples as the initial training samples and sequentially added 10 more samples with the most model uncertainty in each epoch.
As mentioned in Section\ref{sec:closed-loop}, we implement our classifier using MLP model, consisting of 1 input layer, 4 hidden layers, and 1 output layer. ReLU and Adam serve as the activation function and optimizer, respectively.
We also apply the transfer learning technique using a model obtained for the time window [14,16] with historical data.  
In transfer learning mode, the first hidden layer is frozen while the remaining layers are kept trainable.
We compare the performance of the proposed work with the benchmark using a random sampling approach with the same initial model and hyperparameters. The performance of the proposed method without the transfer learning is shown in Fig.~~\ref{fig:case1}(b). Compared to the benchmark shown in Fig.~\ref{fig:case1}(a), the proposed method shows much superior performance, as it pinpoints the boundary of the DSFS much faster and more accurately due to the well-positioned samples.
As shown in Fig.~\ref{fig:case1}(c), with the transferred model, the classifier achieves even better results, despite the fact the DSFS of time [14,16]  (similar to the boundary characterized in epoch 5) is quite visually different. 
In comparison, it can be observed that existing methods~\cite{chen2021leveraging,cui2020network,taheri2022data} that use hyperbox or ellipsoid for inner approximation may be more conservative than our results, as these sets do not fully capture the geometry of the DSFS.
Note that classifying a batch of 1000 samples in GPU with the proposed work on a laptop with Intel(R) UHD Graphics 620 and Core i5-8350U takes only 0.001s. 
Meanwhile, checking one sample using the traditional simulation-based method with Mosek 9.1.9 takes about 0.2s. 
The performance of our approach is credited to the simple operations utilized in the MLP model. During prediction, the computations primarily consist of basic matrix operations and element-wise manipulations. These can be effectively parallelized across multiple processing units, such as GPUs.



\begin{figure}[t]
\centering
\includegraphics[width=0.85\linewidth]{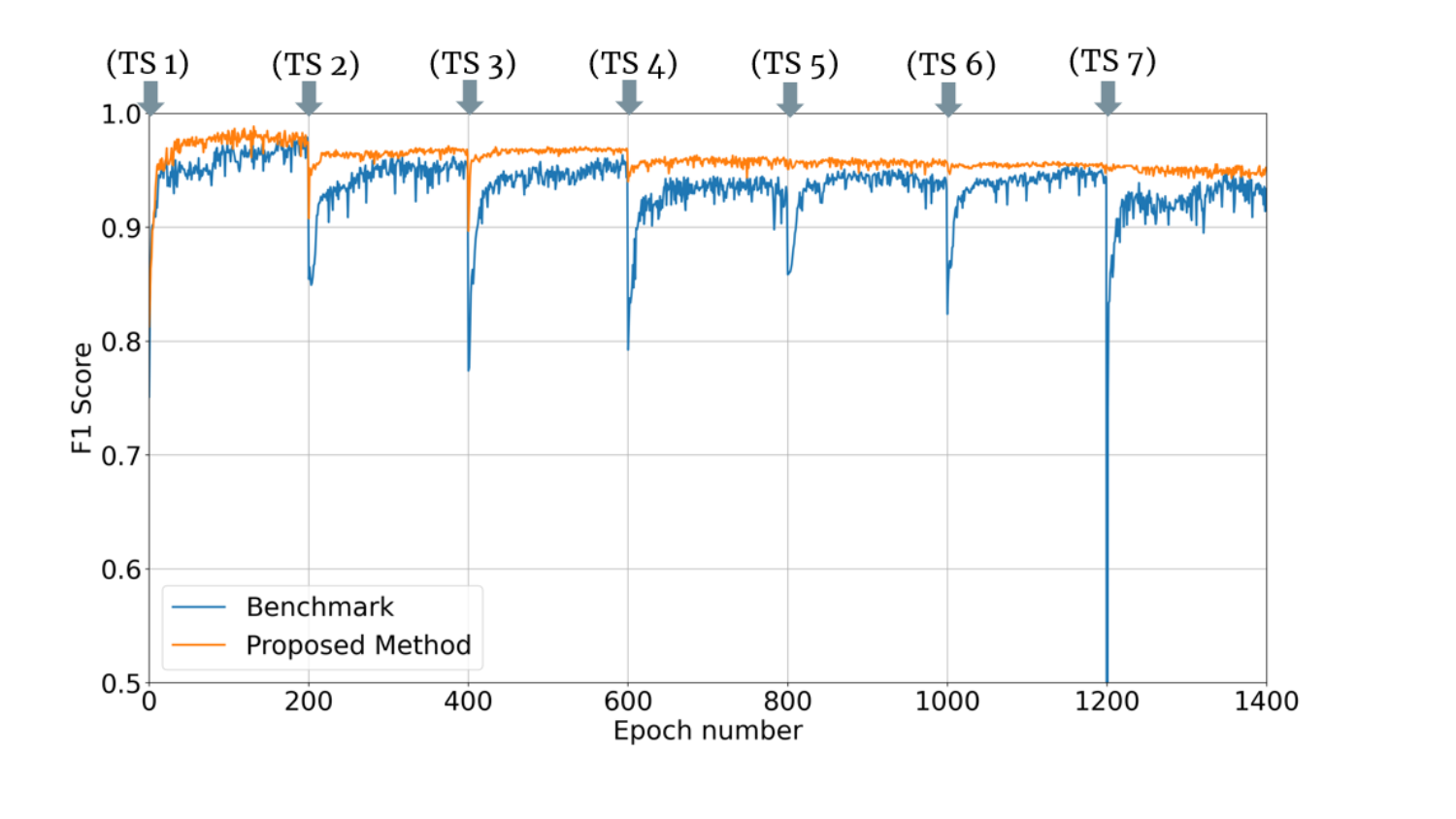}
\caption{Rolling-horizon DSFS estimation results.}
\vspace{-0.15in}
\label{fig:case2}
\end{figure}

We then show the scalability of the proposed work. 
We consider such a scenario that the distribution system estimates the flexibility four-time steps ahead in a rolling horizon, from hour 8 to hour 14. 
From TS 2, we initialize the classifier model with the one obtained from the previous TS. From Fig.~\ref{fig:case2}, the benchmark with the random sampling approach can only achieve the same level of accuracy as ours with almost 10 times more training iterations in TS 1, and cannot keep up for all the following TSs anymore.  

To study the adaptability of our model against noise, we consider the DER injection uncertainty.
We introduce varying levels of uncertainty into the PV system and loads on each node, generating 1000 samples for each uncertainty level as a new test dataset. 
Table~\ref{tab:noise_performance} shows the F1 score performance of our classifier across a range of DER injection uncertainty levels, spanning from 3\% to 40\%.
It can be observed that even at an uncertainty level of 20\%, our model consistently achieves an F1 score exceeding 0.95, indicating its robustness. Moreover, at a heightened uncertainty level of 40\%, the F1 score remains high at 0.88. 
The results show the notable adaptability of our model against uncertainty. The observation of the decreasing estimation accuracy also implies that uncertainties indeed affect the geometry of the DSFS.

\begin{table}[]
\caption{F1 Score of Classifier under Different DER Injection Uncertainty Levels, Ranging From 3\% to 40\%}
\label{tab:noise_performance}
\resizebox{\linewidth}{!}{%
\begin{tabular}{cccccc}
\hline
Uncertainty Level               & 3\%   & 10\%  & 20\%  & 30\%  & 40\%  \\ \hline
F1 score & 0.990 & 0.981 & 0.957 & 0.930 & 0.882 \\ \hline
\end{tabular}%
}
\end{table}

\vspace{-0.1in}

\vspace{-0.05in}
\section{Conclusion}
We propose a data-driven approach to approximate the DSFS. It involves using a new network-informed method to train a classifier that only needs to use scalable matrix operations for the approximation. 
We propose a numerically efficient training strategy that uses the network information and the accumulated knowledge about the sample space to accelerate the training. 
Case studies based on the SCE system verify the validity and value of the proposed work.

\vspace{-0.1in}

\bibliographystyle{IEEEtran}
\bibliography{distri_flex.bib}
\end{document}